\documentclass[11pt]{article}

\usepackage{amsmath,amsthm,amsfonts,amssymb,amscd}

\usepackage{graphicx}

\newtheorem{theoremN}{Theorem}[section]
\newtheorem{propositionN}[theoremN]{Proposition}
\newtheorem{corollaryN}[theoremN]{Corollary}
\newtheorem{lemmaN}[theoremN]{Lemma}

\newcommand{\eop}{\hfill $\Box$}

\newcommand{\R}{\mathbf R}
\newcommand{\Q}{\mathbf Q}

\newcommand{\E}{\mathbf E}

\begin{document}

\author{Paulo R. C. Ruffino\thanks{%
ACKNOWLEDGEMENT: Part of the results in this paper appears firstly in the
author's Ph.D. thesis at the University of Warwick. I want to express my
gratitude to my supervisor Prof. David Elworthy for the endless help and
advice. Research supported by CNPq grant $n^{\circ }$ $300670/95$-$8$. } \\
Departamento de Matem\'{a}tica,\\
Universidade Estadual de Campinas\\
Cx. Postal 6065\\
13.081-970 Campinas-SP, Brazil}
\title{A Sampling Theorem for Rotation Numbers of Linear Processes in 
${\R}^{2}$}

\maketitle

\begin{abstract}
We prove an ergodic theorem for the rotation number of the composition of a
sequence os stationary random homeomorphisms in $S^{1}$. In particular, the
concept of rotation number of a matrix $g\in Gl^{+}(2,{\R})$ can be
generalized to a product of a sequence of stationary random matrices in 
$Gl^{+}(2,{\R})$. In this particular case this result provides a
counter-part of the Osseledec's multiplicative ergodic theorem which
guarantees the existence of Lyapunov exponents. A random sampling theorem is
then proved to show that the concept we propose is consistent by
discretization in time with the rotation number of continuous linear
processes on ${\R}^{2}.$
\end{abstract}

\vspace{3mm} \noindent {\em Key words and phrases:} rotation number, product
of random matrices, sampling theorem.

\vspace{2mm} \noindent {\em 1991 AMS classification scheme numbers:} 60H30,
58F11.

\section{Introduction}

The rotation number of an orientation--preserving homeomorphism in the
circle $f:S^{1}\rightarrow S^{1}$ describes the average rotation that $f$
does in $S^{1}$, either if we consider this average in time for the orbit of
a given starting point or the average on the angular displacement with
respect to an $f$--invariant measure on $S^{1}$ ({\it cf.} Corollary \ref
{4-corolario}). Considering the radial projection ${\R}^{2}\setminus
\{0\}\rightarrow S^{1}$, this concept extends naturally to the space of
linear operators on ${\R}^{2}$ with positive determinant $Gl^{+}(2,{\R})$. So, 
on one hand the Lyapunov exponents give the asymptotic exponential
rate of increasing in the radial coordinate, on the other hand, rotation
number completes the information about the long time behaviour of the system
giving the asymptotic angular behaviour.

Let $g$ be a matrix in $Gl^{+}(2,{\R})$, when we write its eigenvalues
as 
\[
\lambda _{1}=e^{a+bi}\ \ \ \ \mbox{and}\ \ \ \ \ \lambda _{2}=e^{a^{\prime
}-bi}\ , 
\]
with $a=a^{\prime }$ if $b\in (0,\pi )$, the Lyapunov spectrum is $%
\{a,a^{\prime }\}$ and the rotation number is 
\[
\rho (g)\ =\ \pm \ \frac{b}{2\pi }\ \ , 
\]
with the sign depending on the orientation chosen in the plane (see
Proposition \ref{6-eigenvalues}). This twined relation between Lyapunov
exponents and rotation numbers suggests the questions:

\begin{description}
\item[i)]  Is there any extension of the concept of rotation number for the
product of a sequence of stationary random matrices in $Gl^{+}(2,{\R})$?

\item[ii)]  If it does, can we assure its existence almost surely like
Oseledec's multiplicative ergodic theorem does for their Lyapunov exponents?

\item[iii)]  Is this concept consistent by discretization in time (sampling)
with the rotation number of continuous linear processes in ${\R}^{2}$ ?
\end{description}

The main purpose of this paper is to prove that the answers to these
questions are affirmative. In section 2 we review the classical definition
of rotation number and we give an alternative interpretation in terms of 
{\em ordered lifted orbit} of a point $p$ in $S^{1}$ (Proposition \ref
{1-lifted orbit}). In section 3 we extend the definition of rotation numbers
for a sequence of orientation--preserving homeomorphisms in $S^{1}$, and
although the main interest in this paper is rotation number for linear
processes, we prove an ergodic theorem of existence of these numbers in a
more general context: for the composition of a sequence of stationary random
orientation--preserving homeomorphisms in $S^{1}$. In section 4 we study the
case of a product of random matrices. A random sampling theorem is then
proved in section 5 to show that the concept we propose is consistent with
the rotation number of continuous linear cocycles in ${\R}^{2}$.

\section{Classical Definition}

Let $Hom^{+}(S^{1})$ denote the set of homeomorphisms which preserve
orientation in the circle $S^{1}$ . Given $f\in Hom^{+}(S^{1})$, we say that
a continuous monotone increasing function $F:{\R}\rightarrow {\R}$
is a lift of $f$ if $f\circ \pi =\pi \circ F$, where $\pi $ is the canonical
projection $\pi :{\R}\rightarrow S^{1}$, $x\longmapsto e^{2\pi xi}$. Two
lifts of the same homeomorphism will differ by an integer translation. Given
a lift $F$ of $f$ then $F(x+1)=F(x)+1$, which implies that $(F-Id)$ is
periodic with period 1 and 
\[
\max_{x\in {\R}}\{(F-Id)(x)\}-\min_{x\in {\R}}\{(F-Id)(x)\}<1\ . 
\]

We set $F^{(n)}=F\circ F\circ \ldots ,F$, $n$--times. The classical
definition of rotation number of $f$ is, then: 
\begin{equation}
\rho (f)=\lim_{n\rightarrow \infty }\frac{F^{(n)}(x)-x}{n}\pmod{1}
\label{eq.1}
\end{equation}
and it is well known that this limit exists and is independent of the
initial point $x$ and of the lift $F$ chosen. The homeomorphism $%
f:S^{1}\rightarrow S^{1}$ has periodic point if and only if $\rho (f)\in
{\Q}$ , and if $\rho (f)\in {\R}\setminus $ ${\Q}$ then $f$ is
semi--conjugate with $R_{\rho (f)}$, the rotation by the angle $\rho (f)$
(see {\it e.g.} Nitecki\cite{Nitecki}).

Given a point $p\in S^{1}$ and $x\in (-1/2,1/2]$ such that $\pi (x)=p$, we
construct inductively the following real, monotone non-decreasing sequence: $%
\theta _{p,0}=x$ and $\theta _{p,n}$ is the smallest real number such that $%
\theta _{p,n}\geq \theta _{p,n-1}$ and $\pi (\theta _{p,n})=f^{(n)}(p)$. The
sequence $\{\theta _{p,n}\}_{n\geq 0}$ is the ordered orbit of $p$ in the
covering space ${\R}$, it will be called the {\em ordered lifted orbit
of $p$}. Note that $0\leq \theta _{p,n}-\theta _{p,n-1}<1$ for all positive
integer $n$.

\begin{propositionN}
\label{1-lifted orbit}\label{prop2.1}The rotation number $\rho (f)$ is the
asymptotic average rate of increase of the ordered lifted orbit of $p$,
independently of $p$ chosen, {\it i.e.}, 
\begin{equation}
\rho (f)=\lim_{n\rightarrow \infty }\frac{\theta _{p,n}-\theta _{p,0}}{n}%
\pmod{1}.  \label{eq.2}
\end{equation}
\end{propositionN}

\noindent {\bf Proof:} We assume firstly that $f$ has no fixed point in $%
S^{1}$. Without loss of generality, take the unique lift $F$ such that $%
0<(F-Id)(x)<1$ for all $x\in {\R}$, then, since $F$ is strictly
increasing, we have by construction that $\theta _{p,n}=F^{(n)}(x)$ for all $%
n\geq 0$, which implies equation (2).

Suppose now that $f$ has at least one fixed point in $S^{1}$, hence $\rho
(f)=0$. Take the unique lift $F$ which has fixed point in ${\R}$. If $%
F(x)\geq x$, then again $\theta _{p,n}=F^{(n)}(x)$ and this sequence
converges to the fixed point of $F$ which is the nearest from $x$ by above.
So, the limit in equation (2) vanishes and the equality holds. Finally, if $%
F(x)<x$, let $y$ be the fixed point of $F$ which is the nearest from $x$ by
below, then $F^{(n)}(x)$ is a strictly decreasing sequence converging to $y$
and $\theta _{p,1}=F(x)+1,\ldots ,\theta _{p,n}=F(\theta _{p,n-1})+1$, hence

\[
\lim_{n\rightarrow \infty} \frac{\theta_{p,n}-\theta_{p,0}}{n} = \lim_{n
\rightarrow \infty} \frac{F^{(n)}(x)+n-x}{n}=0 \pmod{1}. 
\]
Moreover, since the rotation number is independent of the initial point $x$,
then the limit in equation ($\ref{eq.2}$) is also independent of $p$.

\hfill $\Box$

\noindent {\sc Remark:} If the homeomorphisms $f$, $g:S^{1}\rightarrow S^{1}$
commute then $\rho (f\circ g)=\rho (f)+\rho (g)\pmod{1}$ .

\section{Rotation Number for Composition of Random Homeomorphisms in the
Circle}

Now consider a sequence of orientation--preserving homeomorphisms $f_{1},\
f_{2},\ldots :S^{1}\rightarrow S^{1}$. What we have seen in the last section
suggests to us two approaches to extend the concept of rotation number for
the composition of these homeomorphisms $(f_{n}\circ f_{n-1}\circ \ldots
f_{1})_{n\geq 1}$: \vspace{5mm}\newline
{\sc The first approach. } To consider the composition of a sequence of
lifts $f_{1}^{\sim },\ f_{2}^{\sim },\ldots $ of $f_{1},\ f_{2},\ldots $,
respectively. For the same reason that in the definition of rotation number
for a single homeomorphism $f$ with a lift $F$, when we take the composition 
$F^{(n)}=F\circ \ldots \circ F$ we can not make this composition with
distinct lifts, say $F_{1}\circ F_{2}\ldots \circ F_{n}$ with $%
(F_{j})_{j\geq 1}$ a sequence of distinct lifts of $f$, here, in order to
have the same kind of compatibility among the lifts $f_{1}^{\sim
},f_{2}^{\sim },\ldots $, we will state that all of them start in the same
interval, say $f_{j}^{\sim }(0)\in (-1/2,1/2]$ for all $j$ positive integer.
Set $F_{n}^{\sim }=f_{n}^{\sim }\circ \ldots f_{1}^{\sim }$ and given $x\in 
{\R}$, define:

\begin{equation}
\rho (f_{1},f_{2},\ldots )=\lim_{n\rightarrow \infty }\frac{F_{n}^{\sim
}(x)-x}{n}\pmod{1},  \label{eq.4}
\end{equation}
when the limit exists. \vspace{5mm}\newline
{\sc The second approach.} Define the rotation number of the composition $%
f_{n}\circ \ldots \circ f_{1}:S^{1}\rightarrow S^{1}$, $n=1,2,\ldots $ via
the ordered lifted orbit of a point $p\in S^{1}$. Namely if $x\in (-1/2,1/2]$
such that $\pi (x)=p$ then construct inductively the sequence: $\theta
_{p,0}=x$, and $\theta _{p,n}$ is the smallest real number such that $\theta
_{p,n-1}\leq \theta _{p,n}$ and $\pi (\theta _{p,n})=f_{n}\circ \ldots \circ
f_{1}(p)$, and define 
\begin{equation}
rot(f_{1},f_{2},\ldots )(p)=\lim_{n\rightarrow \infty }\frac{\theta
_{p,n}-\theta _{p,0}}{n}\pmod{1}  \label{eq.3}
\end{equation}
when the limit exists. \vspace{5mm}\newline
\indent Next proposition shows that the first approach is particularly
interesting for our purposes.

\begin{propositionN}
\label{2-invariancia}\label{prop3.1}If the rotation number $\rho
(f_{1},f_{2},\ldots )$ of equation ($\ref{eq.4}$) exists then it is
independent of the initial point $x\in {\R}$.
\end{propositionN}

\noindent {\bf Proof:} The proof goes similarly to the proof for the single
homeomorphism case. Arguing inductively, $F_{n}^{\sim }$ as defined above is
a lift of $f_{n}\circ \ldots \circ f_{1}$, hence, for all positive integer $%
n $: 
\[
\max_{x\in {\R}}\{(F_{n}^{\sim }-Id)(x)\}-\min_{x\in {\R}%
}\{(F_{n}^{\sim }-Id)(x)\}<1\ , 
\]
then, 
\[
|F_{n}^{\sim }x-F_{n}^{\sim }y|\leq |(F_{n}^{\sim }x-x)-(F_{n}^{\sim
}y-y)|+|x-y|\leq |x-y|+1, 
\]
therefore 
\[
\lim_{n\rightarrow \infty }\left( \frac{F_{n}^{\sim }x}{n}-\frac{F_{n}^{\sim
}y}{n}\right) =0; 
\]
hence if the limit of equation ($\ref{eq.4}$) exists, it is independent of $%
x $.

\hfill $\Box$

\noindent {\sc Example 1.} For simplicity we parametrize $S^1$ by $x\mapsto
e^{2\pi x i}$ with $x\in (-1/2, 1/2]$. Let $f_1,f_2,f_3 $ and $f_4$ be such
that $f_i(0)=0$, i=1,2,3,4, and $f_1(1/8)=3/8$, $f_2(3/8)=-3/8$, $%
f_3(-3/8)=-1/8$ and $f_4(-1/8)=1/8$. Consider the sequence $f_1,f_2, \ldots$
such that $f_n=f_{(n\ \bmod\ 4 )}$. For the first approach $\rho
(f_1,f_2,\ldots )=0$ since $0$ is a fixed point of $f_n\circ \ldots \circ
f_1 $. For the second approach $rot(f_1,f_2, \ldots )(0)=0$, nevertheless,
if $p=1/8$, then 
\[
\theta_{p,n}=\frac{2n+1}{8}\ , 
\]
which yields $rot(f_1,f_2,\ldots )(p)=1/4$.

This example shows not only that the two approaches lead to different
numbers but also that in the second one this number depends on the starting
point $p$. From now on, whenever we refer to rotation number, we will mean
the first approach. \vspace{2mm}\newline
\noindent {\sc Remark:} Proposition \ref{1-lifted orbit} shows that for a
single homeomorphism $f$, the two approaches coincide. \vspace{1mm}\newline

Once stated the compatibility condition that any lift $F$ will be taken with 
$F(0)\in (-1/2,1/2]$, we can set a metric in $Hom^{+}(S^{1})$ given by the
uniform metric in the lifts. This definition creates natural discontinuities
at homeomorphisms $f$ with lift $F(0)=1/2$.

Although in this paper our main interest is rotation numbers for linear
systems, an ergodic theorem of existence of such numbers can be stated in a
more general context: for composition of a sequence of stationary random
homeomorphisms in $Hom^{+}(S^{1})$.

To introduce formally the set up, consider a probability space ($\Omega,%
{\cal F}, {\rm I\!P} $) and $\theta:\Omega \rightarrow \Omega$ a
measure--preserving transformation on $\Omega$. Let us assume throughout for
simplicity that $\theta$ is ergodic. Let $f: \Omega \rightarrow Hom^+(S^1)$
be a random homeomorphism in $S^1$, we shall consider the sequence of
stationary random homeomorphisms given by $f_n=f \circ \theta^{n-1}$. The
main result of this section guarantees the existence ${\rm I\!P}$-almost
surely of the rotation number $\rho (f, \theta) $ of the discrete random
dynamical system on $S^1$ over $\theta$ given by the composition $(f_n\circ
f_{n-1}\circ \ldots f_1)_{n\geq 1}$.

For each $n\geq 1$, we will write the lift $f_{n}^{\sim }$ as 
\begin{equation}
f_{n}^{\sim }(\omega ,x)=x+\delta _{n}(\omega ,x)\ ,  \label{eq:psi}
\end{equation}
where the function $\delta _{n}(\omega ,x)$ is periodic in the variable $x$
with period 1 and $|\delta _{n}(\omega ,x)|<3/2$, for all $x$. In the sequel
it will be convenient to write the $\delta _{n}$'s as functions on $S^{1}$,
so for each $n$ we define $\beta _{n}:\Omega \times S^{1}\rightarrow
(-3/2,3/2)$ by 
\begin{equation}
\beta _{n}(\omega ,\pi (x))=\delta _{n}(\omega ,x)\ .  \label{eq:beta}
\end{equation}

We denote by $\mu$ an invariant probability measure on $\Omega\times S^1$
for the skew product map $\Theta(\omega,s)=(\theta(\omega), f(\omega, s))$.
The invariant measure $\mu$ factorizes as $\mu(ds, d\omega) =
\nu_{\omega}(ds){\rm I\!P}(d\omega)$ (see, {\it e.g.} Crauel \cite{Crauel}
or Arnold \cite{Arnold}).

With this set up we state the following ergodic theorem:

\begin{theoremN}
\label{3-teo ergodico}\label{thm3.2}Consider the discrete random dynamical
system on $S^{1}$ over $\theta $ given by $(f_{n}\circ f_{n-1}\circ \ldots
f_{1})_{n\geq 1}$ where $f_{n}=f\circ \theta ^{n-1}$. Then, its rotation
number $\rho (f,\theta )$ exists ${\rm I\!P}$-a.s. and satisfies: 
\begin{equation}
\rho (f,\theta )={\rm I\!E}\left[ \ \int_{S^{1}}\beta _{1}(\omega ,s)\ d\nu
_{\omega }(s)\ \right] \ \pmod{1}\ \ \ \ a.s.  \label{D.RN.station}
\end{equation}
where $\mu =\nu _{\omega }(ds){\rm I\!P}(d\omega )$ is a (not necessarily
ergodic) invariant probability measure on $\Omega \times S^{1}$ for the skew
product map $\Theta $.
\end{theoremN}

\noindent {\bf Proof:} By construction: 
\[
\beta _{i}\Bigl(\omega ,f_{i-1}\circ \ldots f_{1}(p)\Bigr)=\beta _{1}\circ
\Theta ^{i-1}\Bigl(\omega ,p\Bigr)\ , 
\]
for $i\geq 1$. Assume that $x\in {\R}$ is an initial point with $\pi
(x)=p$, then by equation (\ref{eq:psi}) and induction on $n$: 
\[
f_{n}^{\sim }\circ \ldots \circ f_{1}^{\sim }(\omega ,x)=x+\delta
_{1}(\omega ,x)+\delta _{2}\Bigl(\omega ,x+\delta _{1}(\omega
,x)\Bigr)+\ldots 
\]
\[
+\delta _{n}\Bigl(\omega ,x+\delta _{1}(\omega ,x)+\delta _{2}(\omega
,x+\delta _{1}(\omega ,x))+\ldots +\delta _{n-1}(\ldots )\Bigr)\ . 
\]
Moreover, for any $i=1,\ldots ,n$, by periodicity of $\delta _{i}$ (equation 
\ref{eq:beta}): 
\begin{eqnarray}
\delta _{i}\Bigl(\omega ,x+\delta _{1}(\omega ,x)+\ldots +\delta
_{i-1}(\ldots )\Bigr) &=&\delta _{i}\Bigl(\omega ,f_{i-1}^{\sim }\circ
\ldots \circ f_{1}^{\sim }(\omega ,x)\Bigr)  \nonumber \\
&=&\beta _{i}\Bigl(\omega ,f_{i-1}\circ \ldots \circ f_{1}(\omega ,p)\Bigr) 
\nonumber \\
&=&\beta _{1}\circ \Theta ^{i-1}\Bigl(\omega ,p\Bigr)\ ,  \label{eq:delta}
\end{eqnarray}
hence: 
\[
f_{n}^{\sim }\circ \ldots f_{1}^{\sim }(\omega ,x)=x+\sum_{i=1}^{n}\beta
_{1}\circ \Theta ^{i-1}(\omega ,p)\ . 
\]
If the measure $\mu $ is ergodic then by the Birkhoff's ergodic theorem: 
\[
\lim_{n\rightarrow \infty }\ \frac{1}{n}f_{n}^{\sim }\circ \ldots
f_{1}^{\sim }(x)={\rm I\!E}\left[ \ \int_{S^{1}}\beta _{1}(\omega ,s)\ d\nu
_{\omega }(s)\ \right] \ \ \ \ \ \mu -a.s.. 
\]
Note that once $\theta $ is ergodic on $\Omega $, last formula says that
there exists a subset $\Omega ^{\prime }\subset \Omega $ of probability one
such that for each $\omega ^{\prime }\in \Omega ^{\prime }$, there exists $%
p\in S^{1}$ with the equality above holding for $(\omega ^{\prime },p)$. The
fact that the rotation number is independent of the initial point in $S^{1}$
(Proposition \ref{2-invariancia}) implies that the equality also holds for $%
(\omega ^{\prime },s)$ for any $s\in S^{1}$. Hence, this dynamical
restriction implies that the integral above does not depend on the ergodic
measure chosen and can be taken with respect to any invariant probability
measure, once they are convex combinations of the ergodic measures.

\hfill $\Box$

In particular, if the sequence $\{ f_n\}_{n\geq 1}$ is i.i.d. then the
process $(f_n\circ f_{n-1}\circ \ldots f_1(s))_{n\geq 1}$ in the circle $S^1$
is Markovian. In this case an invariant measure $\mu$ is given by the
product measure $\nu\otimes{\rm I\!P}$ where $\nu$ is a stationary
probability measure for the Markov process in $S^1$. Hence, in this case the
rotation number is given by: 
\begin{equation}
\rho(f ) = \int_{S^1} {\rm I\!E} \ [\ \beta_1(\omega,s)\ ]\ d\nu (s) \ %
\pmod{1} \ \ \ \ a.s..  \label{D.RN.iid}
\end{equation}

We finish this section with a corollary whose proof follows naturally in
this context, although it does not require a sequence of homeomorphisms. It
shows that the rotation number of a homeomorphism $f$ on $S^1$ is the
average of angular displacement of the points in the circle with respect to
any probability measure which is preserved by $f$. The proof comes directly
from formula ($\ref{D.RN.iid}$).

\begin{corollaryN}
\label{4-corolario}\label{cor3.3}If $f\in Hom^{+}(S^{1})$ preserves the
probability measure $\nu $ on $S^{1}$ then 
\[
\rho (f)=\int_{S^{1}}\beta (s)\ d\nu (s)\ \ \ \pmod{1}\ ,
\]
where $\beta :S^{1}\rightarrow {\R}$ is a continuous function on $S^{1}$
which gives a lift $F(x)=x+\beta \left( \pi (x)\right) \ $ for $f$.
\end{corollaryN}

\section{Product of Random Matrices in $Gl^{+}(2,{\R})$}

We will denote by $\psi _{g}:S^{1}\rightarrow S^{1}$ the action of the
matrix $g\in Gl^{+}(2,{\R})$ over the circle $S^{1}$, {\it i.e.} $\psi
_{g}(x)=\frac{gx}{\Vert gx\Vert }$. The rotation number of $\psi _{g}$ will
be denoted simply by $\rho (g)$.

The next proposition shows that the rotation number of a matrix $g\in
Gl^{+}(2,{\R})$ is a concept twined with its Lyapunov numbers in the
sense that if $\lambda _{1},\lambda _{2}\in {\bf C}$ are the eigenvalues of $%
g$ then the logarithm of their modulus $|\lambda _{1}|$ and $\ |\lambda
_{2}| $ give the Lyapunov exponents and their arguments $\pm \arg (\lambda
_{1})$ give the rotation number. To prove this proposition we use the
following lemma which states the invariance of the rotation number by
conjugacy:

\begin{lemmaN}
\label{lemma4.1}If $f,h\in Hom^{+}(S^{1})$ then $\rho (h\circ f\circ
h^{-1})=\rho (f)$.
\end{lemmaN}

\noindent {\bf Proof:} Let $H$ be a lift for $h$ and consider its inverse $%
H^{-1}:{\R}\rightarrow {\R}$. Applying $h^{-1}$ in both sides of $%
h\circ \pi =\pi \circ H$ we check that $H^{-1}$ is also a lift for $h^{-1}$.
For a lift $F$ of $f$ we have $(h\circ f\circ h^{-1})\circ \pi =\pi \circ
(H\circ F\circ H^{-1})$, {\it i e}, $H\circ F\circ H^{-1}$ is a lift for $%
h\circ f\circ h^{-1}$ as well. The result follows immediately from
definition: \label{5-lemma} 
\[
\lim_{n\rightarrow \infty }\frac{(H\circ F\circ H^{-1})^{n}(x)-x}{n}%
=\lim_{n\rightarrow \infty }\frac{H\circ F^{(n)}\circ H^{-1}(x)-x}{n}=\rho
(f)\pmod{1}\ 
\]
since $H(x)=Id(x)+\delta (x)$ where $\delta $ is a bounded periodic function.

\hfill $\Box$

\noindent{\sc Remark:} If we allow $h$ to be a reverse--orientation
homeomorphism in $S^1$ then its lift is a continuous monotone strictly
decreasing functions, and it can be written as $H(x)= -Id(x)+\delta(x)$, for
some periodic bounded function $\delta$. So, in this case, the effect in the
rotation number of the conjugacy by $h$ is the change of the sign: 
\begin{equation}
\rho(h\circ f\circ h^{-1})= -\rho(f) \pmod{1}\ .  \label{eq.inv.signal}
\end{equation}
It corresponds of considering the rotation in the opposite (clockwise)
direction. \vspace{3mm}

Note that actually it does not matter in which interval we consider the
rotation number $\rho (f)$ since it represents an equivalence class for the $%
\pmod{1}$ relation. However expressions like equation ($\ref{eq.inv.signal}$%
) would look more natural if we write $\rho (f)$ in the symmetric interval $%
(-1/2,1/2]$. Another advantage of the representation of the equivalence
classes in this interval is that it becomes easily comparable with rotation
number for continuous process in $S^{1}$ via discretization, provided we
take samples of the process at intervals of time $T>0$ small enough (see
Theorems \ref{7-sampling thm} and \ref{8-determ.sampling}).

\begin{propositionN}
\label{6-eigenvalues}\label{prop4.2}Let $g$ be a matrix in $Gl^{+}(2,{\R}%
)$ and $\lambda _{1},\ \lambda _{2}\in {\bf C}$ be its eigenvalues. Then 
\[
\rho (g)=\pm \frac{1}{2\pi }\arg (\lambda _{1})\ \ \ ,
\]
where $\arg (\lambda _{1})\in (-\pi ,\pi ]$ is the argument of the complex
number $\lambda _{1}$, and the sign will depend on the orientation chosen in 
${\R}^{2}$.
\end{propositionN}

\noindent {\bf Proof:} The proof comes immediately from Lemma 4.1: write $g$
in its Jordan decomposition form $P\Lambda P^{-1}$, with $P\in Gl^{+}(2,%
{\R}^{2})$; then $\rho (g)=\rho (\Lambda )$.

\hfill $\Box$

Considering the action of the group $Gl^{+}(2,{\R})$ on the circle $%
S^{1} $, the concept of rotation number for the composition of a sequence of
homeomorphisms is extended to the product of matrices in $Gl^{+}(2,{\R})$%
. By Theorem \ref{3-teo ergodico}, for a stationary sequence of random
matrices (in particular for i.i.d.) this number exists ${\rm I\!P}$-almost
surely, in the same way that the existence of the Lyapunov exponents in this
case is assured by the Oseledec's multiplicative ergodic theorem (see {\it %
e.g.} Ruelle \cite{R1} or Furstenberg and Kifer \cite{FurstK} for a
non--random filtration approach).

To set up the notation, let $Y:\Omega \rightarrow Gl^{+}(2,{\R})$ be a
random matrix and $\theta :\Omega \rightarrow \Omega $ an ergodic
transformation on $\Omega $. We shall consider the sequence of stationary
random matrices given by $Y_{n}=Y\circ \theta ^{n-1}$. Let $\Psi _{Y_{n}}:%
{\R}\rightarrow {\R}$ be the lift of $\psi _{Y_{n}}\in
Hom^{+}(S^{1}) $, for each positive integer $n$. As before, we write 
\[
\Psi _{Y_{n}}(x)=x+\delta _{n}(\omega ,x)\ \ . 
\]
Note that here, because of the linearity of $Y_{n}$, the functions $\delta
_{n}(\omega ,x)$ has period $1/2$, also any invariant measure $\mu =\nu
_{\omega }(ds){\rm I\!P}(d\omega )$ on $\Omega \times S^{1}$ can be
considered such that $\nu _{\omega }$ is a measure on the projective space $P%
{\R}^{1}$. We write $\beta _{n}(\omega ,\pi (x))=\delta _{n}(\omega ,x)$%
, then by Theorem \ref{3-teo ergodico} the rotation number of the product of
this sequence of random matrices exists and satisfies: 
\begin{equation}
\rho (Y,\theta )={\rm I\!E}\left[ \ \int_{S^{1}}\beta _{1}(\omega ,s)\ d\nu
_{\omega }(s)\ \right] \ \ \ \ \pmod{1}\ \ \ \ \ a.s..
\label{RN.prm.station}
\end{equation}
In the i.i.d case ({\it cf.} formula ($\ref{D.RN.iid}$)): 
\begin{equation}
\rho (Y)=\int_{S^{1}}{\rm I\!E}\ [\ \beta _{1}(\omega ,s)\ ]\ d\nu (s)\ \ \
\ \ \pmod{1}\ \ \ \ \ a.s.,  \label{RN.prm.iid}
\end{equation}
with $\nu $ a stationary probability measure on $S^{1}$ for the Markov
process.

We present two simple but illustrative examples:

\noindent {\sc Example 2.} Let $Y_{1},Y_{2},\ldots $ be a sequence of i.i.d.
random matrices such that the support of the common distribution is a subset
of $SO(2,{\R})$. Then for each $n$, we can associate a real random
variable $\lambda _{n}(\omega )\in (-1/2,1/2]$ such that $Y_{n}(\omega )$ is
the rotation by $2\pi \lambda _{n}(\omega )$, with $\lambda _{n}(\omega )$
i.i.d. Hence the lift $\Psi $ satisfies $\Psi _{Y_{n}(\omega )}(x)=x+\lambda
_{n}(\omega )$, for all $n$, and by the law of large numbers:

\[
\rho(Y_1, Y_2, \ldots )= \lim_{n\rightarrow \infty } \ \frac {1}{n} \
\sum_{j=1}^n \ \lambda_j(\omega ) = \ {\rm I\!E} \ [\lambda_1 ]\ \ \ a.s. . 
\]
On the other hand, since $\beta_1(\omega,x) \equiv \lambda_1(\omega)$ for
all $s\in S^1$ and measure $\nu$ is the normalized Lebesgue measure on $S^1$%
, one easily verifies formula ($\ref{RN.prm.iid}$).

\hfill $\Box$

\noindent {\sc Example 3.} Let $Y_{1},Y_{2},\ldots $ be a sequence of i.i.d.
random upper triangular matrices in $Gl^{+}(2,{\R})$. Set $\
a_{n}(\omega )$ for the 1,1--entry of $Y_{n}(\omega )$, then $a_{n}\neq 0$
and the sequence $(a_{n})_{n\geq 1}$ is i.i.d.. By induction on $k$ the
lifts $\Psi _{Y_{n}}$ are such that, for $k$ positive integer: 
\[
\Psi _{Y_{n}}(k)=\left\{ 
\begin{array}{ll}
k\ \ \ \  & \ \ \ \ \mbox{ if }\ a_{n}>0\ \  \\ 
k+1/2\ \ \ \  & \ \ \ \ \mbox{ if }\ a_{n}<0\ \ ,
\end{array}
\right. 
\]
and 
\[
\Psi _{Y_{n}}(k+1/2)=\left\{ 
\begin{array}{ll}
k+1/2\ \ \ \  & \ \ \ \ \mbox{ if }\ a_{n}>0\ \  \\ 
k+1\ \ \ \  & \ \ \ \ \mbox{ if }\ a_{n}<0\ \ .
\end{array}
\right. 
\]
So, $\Psi _{Y_{n}}(k)$ ``rotates'' half of the circle $S^{1}$ ($\pi \circ
\Psi _{Y_{n}}(k)$ goes to its antipode) if $a_{n}<0$ and ``does not rotate''
(fixed) if $a_{n}>0$. So one easily calculates the rotation number starting
at $x=0$:

\begin{eqnarray}
\rho (Y_{1},Y_{2},\ldots )=\lim_{n\rightarrow \infty }\frac{\Psi
_{Y_{n}}\circ \ldots \circ \Psi _{Y_{1}}(0)}{n} &=&\lim_{n\rightarrow \infty
}\frac{1}{n}\ \sum_{j=1}^{n}\ \frac{1}{2}\mbox{{\bf 1}}_{\{a_{j}<0\}} 
\nonumber \\
&&  \nonumber \\
&=&\frac{1}{2}{\rm I\!P}\left[ a_{1}<0\right] \ \ \ a.s.  \nonumber
\end{eqnarray}
On the other hand, if $\{e_{1},e_{2}\}$ is the standard basis in ${\R}%
^{2}$ then: 
\[
\beta _{1}(\omega ,e_{1})=\beta _{1}(\omega ,-e_{1})=\frac{1}{2}%
\mbox{{\bf
1}}_{\{a_{1}<0\}}\ \ , 
\]
and if $\delta _{(\cdot )}$ denotes the Dirac measure then: 
\[
\nu =\frac{1}{2}\delta _{e_{1}}+\frac{1}{2}\delta _{-e_{1}} 
\]
is an invariant measure on $S^{1}$ (not necessarily ergodic) and one easily
verifies formula ($\ref{RN.prm.iid}$).

\nopagebreak

\hfill $\Box$

\noindent {\sc Remark:}{\em (Rotation number for diffeomorphisms)} Besides
the sampling theorem for continuous linear processes in $R^{2}$ (next
section), another interesting application of the concept of rotation number
for a sequence of random matrices is in non--linear dynamical systems. Let $%
f $ be an orientation preserving diffeomorphism on ${\R}^{2}$ with a
finite invariant measure $\mu $ (or in general in a 2--dimensional manifold
with the support of the invariant measure contained in a neighbourhood where
the tangent bundle $TM$ is parallelizable). The rotation number of $f$ at a
point $p\in {\R}^{2}$, defined by the rotation number of the product of
the sequence of the differential maps 
\[
df^{(n)}(p)=df\Bigl(f^{(n-1)}(p)\Bigr)\circ \ldots \circ df(p)\ , 
\]
$n\geq 1$, gives the average rotation of the directions of the stable (and
unstable) submanifold along the orbit of the point $p$. By formula ($\ref
{RN.prm.station}$) this number exists for $\mu $--almost every point $p$,
moreover it is constant in each ergodic component. It is not our purpose in
this paper to go further in this non--linear analysis; it will be dealt with
elsewhere.

\section{Random Sampling Theorem}

In this section we shall consider the discretization in time of a continuous
linear cocycle on ${\R}^{2}$, in particular, of the solution of a linear
stochastic differential equation. To set up the notation let $(\theta
_{t})_{t\geq 0}$ be a flow of ergodic transformations on the probability
space $(\Omega ,{\cal F},{\rm I\!P})$. A continuous linear (perfect) cocycle 
$\varphi (t,\omega )$ on ${\R}^{2}$ over $\theta $ is a map $\varphi
:\Omega \times {\R}_{\geq 0}\rightarrow Gl^{+}(2,{\R})$ such that
for all $\omega \in \Omega $: \newline
(i) $\varphi (0,\omega )=Id_{2\times 2}$;\newline
(ii) $\varphi (t,\omega )$ is continuous on $t$;\newline
(iii) it has the cocycle property: 
\[
\varphi (t+s,\omega )=\varphi (t,\theta _{s}(\omega ))\circ \varphi
(s,\omega ). 
\]
We deal with perfect cocycles once for every crude cocycle there exists a
perfect cocycle such that they are indistinguishable, see L. Arnold and M.
Scheutzow \cite{ArnScheu} or L. Arnold \cite{Arnold}. This cocycle generates
the following random linear system on ${\R}^{2}$: 
\begin{equation}
x_{t}=\varphi (t,\omega )x_{0}.  \label{eq.cocycle}
\end{equation}

We shall denote by $\xi _{t}(\omega )$ the induced cocycle in the unitary
circle $S^{1}$ given by the radial projection of ${\R}^{2}\setminus
\{0\} $ onto $S^{1}$. Associated with the process $s_{t}=\xi _{t}(\omega
)s_{0}$ there is the real continuous angular process $\alpha _{t}(s_{0})$
(parametrized by the initial condition $s_{0}\in S^{1}$) such that $%
s_{t}=\exp \{i\alpha _{t}(s_{0})\}$, {\it i.e.} $\alpha _{t}(s_{0})$ is the
continuous angular component of $x_{t}$. The {\em rotation number} $\rho
(\varphi )$ of this linear system is the average angular velocity of a
solution starting at $x_{0}\in {\R}^{2}\setminus \{0\}$: 
\[
\rho (\varphi )=\lim_{t\rightarrow \infty }\frac{\alpha _{t}(s_{0})}{t}\ , 
\]
when the limit exists, and if it does, it is independent of the initial
point $s_{0}\in S^{1}$.

For a fixed $T>0$, by the cocycle property of $\varphi (t,\omega )$: 
\[
\varphi \Bigl(nT,\omega \Bigr)=\varphi \Bigl(T,\theta _{(n-1)T}(\omega
)\Bigr)\circ \ldots \circ \varphi \Bigl(T,\theta _{T}(\omega )\Bigr)\circ
\varphi \Bigl(T,\omega \Bigr). 
\]
Next theorem will assure that for a quite large class of cocycles the
rotation number for the product of the random matrices $\{\varphi (T,\theta
_{(n-1)T}(\omega ))\}_{n\geq 1}$ agrees with the rotation number $\rho
(\varphi )$ of the continuous system when $T>0$ approaches zero, up to a
factor of scaling of $1/T$. We shall denote by $\mu =\nu _{\omega }(ds){\rm %
I\!P}(d\omega )$ an invariant probability measure in $S^{1}\times \Omega $
for the skew-product flow in this product space induced by $\xi _{t}(\omega )
$.

We enphasize that, although hypothesis (\ref{rot.cocycle}) looks artificial
and difficult to be verified, we prefer to state the theorem in this general
formulation for cocycles and show latter (sub-sections 5.1, 5.2 and 5.3)
that most of interesting dynamical systems (deterministic, real noise and
stochastic respectively) satisfies naturally this hypothesis.

\begin{theoremN}[Random sampling theorem]
\label{7-sampling thm}\label{thm5.1}Consider a linear cocycle $\varphi
(t,\omega )$ such that there exists the rotation number $\rho (\varphi )$ of
the continuous process and 
\begin{equation}
\rho (\varphi )=\lim_{T\rightarrow 0}\ \frac{1}{T}\ {\rm I\!E}\left[
\int_{S^{1}}\Bigl(\alpha _{T}(s)-\alpha _{0}(s)\Bigr)d\nu _{\omega
}(s)\right] \ \ \ \ a.s..  \label{rot.cocycle}
\end{equation}
Then the rotation number $\rho (\varphi (T,\omega ),\theta _{T})$ of the
sequence of random matrices \newline
$\{\varphi (T,\theta _{T(n-1)}(\omega )\ \}_{n\geq 1}$ satisfies 
\[
\lim_{T\rightarrow 0}\ \frac{1}{T}\ \rho (\varphi (T,\omega ),\theta _{T})\
=\ \rho (\varphi )\ \ \ \ \ \ \ a.s..
\]
\end{theoremN}

\noindent {\bf Proof:} Once stated the compatibility condition that the
lifts should satisfy \newline
$\Psi _{\varphi (T,\cdot )}(0)\in (-1/2,1/2]$ (section 3), the map $T\mapsto
\Psi _{\varphi (T,\omega )}$ is continuous with respect to the uniform
metric defined on $Hom^{+}(S^{1})$ up to the stopping time 
\[
\sigma =\inf \{t\geq 0:\xi _{t}(e_{1})=-e_{1}\}\ . 
\]
So, fix $s_{0}\in S^{1}$ with the corresponding initial angle $\alpha
_{0}(s_{0})$ and take a real $T>0$, then by construction: 
\[
\beta _{1}^{T}(\omega ,s_{0})\ 1_{\{T\leq \sigma \}}=\Bigl(\alpha
_{T}(\omega )-\alpha _{0}(\omega )\Bigr)\ 1_{\{T\leq \sigma \}}. 
\]
Hence by the hypothesis and the Lebesgue's convergence theorem: 
\begin{eqnarray*}
\lim_{T\rightarrow 0}\ \frac{1}{T}\ {\rm I\!E}\left[ \int_{S^{1}}\beta
_{1}^{T}(\omega ,s)\ d\nu _{\omega }(s)\right] &=& \\
&& \\
\lim_{T\rightarrow 0}\ \frac{1}{T}\ {\rm I\!E}\left[
\int_{S^{1}}\Bigl(\alpha _{T}(s)-\alpha _{0}(s)\Bigr)\ 1_{\{T\leq \sigma
\}}\ d\nu _{\omega }(s)\right] &=&\ \rho (\varphi )\ \ \ \ \ \ \ a.s.,
\end{eqnarray*}
and the result follows by formula ($\ref{RN.prm.station}$) and the fact that
an invariant measure for the skew--product flow is also invariant for the
discretized system.

\hfill $\Box$

The next three particular cases show that the theorem holds for most of the
interesting linear cocycles on ${\R}^{2}$.

\subsection{The Deterministic Case}

For deterministic systems it holds a more accurate result than Theorem \ref
{7-sampling thm}:

\begin{theoremN}[Deterministic sampling theorem]
\label{8-determ.sampling}\label{thm5.2}Consider the linear system $\dot{x}=Ax
$, with $A$ a $2\times 2$--matrix and let $\rho $ be its rotation number as
a continuous system. If 
\[
T<\frac{1}{2\rho }\ ,
\]
($T<\infty $ if $\rho =0$) then 
\begin{equation}
\frac{1}{T}\ \rho (\varphi _{T})\ =\ \rho \ \ ,  \label{eq.tnd}
\end{equation}
where $\varphi _{T}\in Gl^{+}(2,{\R})$ is the fundamental solution of
the system at time $T$.
\end{theoremN}

\noindent {\sc Remark:} This theorem says that if the sampling frequency is
greater than $2\rho $ then we can retrace exactly the original frequency
(rotation number) of the continuous system. For a given real signal
(function) $s(t),t\in {\R}$, if $f_{0}$ is the maximum frequency of its
Fourier spectrum, it is well known in the engeneering literature that the
whole spectrum, hence the signal, can be retraced if we sample in time this
signal at a frequency greater than $2f_{0}$. This frequency is called the 
{\em Nyquist's rate} (see {\it e.g.} Papoulis \cite{Papoulis} or Oppenheim
and Schafer \cite{OppScha}). If we identify each frequency in the Fourier
spectrum with the corresponding rotation number of a continuous linear
system on ${\R}^{2}$, then heuristically Theorem \ref{8-determ.sampling}
gives an alternative prove of the property of the Nyquist's rate.

\vspace{3mm}

\noindent {\bf Proof:} Let $\lambda _{1}$, $\lambda _{2}$ be the eigenvalues
of the matrix of coefficients $A$. If $\lambda _{1},\ \lambda _{2}\in {\R%
}$ then there exist fixed points for the continuous flows in $S^{1}$,
therefore $\rho =0$. In this case the eigenvalues of $\varphi _{T}$ are $%
e^{T\lambda _{1}}$ and $e^{T\lambda _{2}}$, hence by Proposition \ref
{6-eigenvalues} $\ \rho (\varphi _{T})=0$ for all $T\in {\R}$.

Now assume that $\lambda _{1},\ \lambda _{2}=a\pm ib$, with $b\neq 0$. It is
well known that in this case $\rho =\pm b$, with the sign depending on the
orientation (see {\it e.g.} Arnold and San Martin \cite{ASM} or San Martin 
\cite{SM}). Then, for an arbitrary $t\geq 0$, the eigenvalues of $\varphi
_{t}$ are $e^{t\lambda _{1}},\ e^{t\lambda _{2}}$ and by Proposition \ref
{6-eigenvalues} again: 
\[
\rho (\varphi _{t})=\pm tb\pmod{1}\ . 
\]
So, for a fixed $T<1/2\rho \ $ we have $\rho (\varphi _{T})=\pm Tb\in
(-1/2,\ 1/2]$, hence 
\[
\frac{1}{T}\ \rho (\varphi _{T})\ =\ \pm \rho \ . 
\]
Given a fixed orientation on $S^{1}$, say anti--clockwise, the agreement in
the sign comes naturally using the ordered lifted orbit (Proposition \ref
{1-lifted orbit}). If $\rho >0$, say, then the lifted orbit for $\psi
_{\varphi _{T}}$ satisfies: 
\[
0<\theta _{x,n}-\theta _{x,n-1}<\frac{1}{2} 
\]
for all $x\in {\R}$ and $n\in {\Bbb N}$, which implies that $\rho
(\varphi _{T})$ is also positive.

\hfill $\Box$

\noindent {\sc Remark:} We emphasize that one of the advantages of fixing $%
\rho (\varphi _{T})$ in the interval $(-1/2,\ 1/2]$ is that it makes sense
to talk about ``positive'' or ``negative'' rotation number. So, the ``change
of the sign'' when we change orientation looks more natural, besides,
equality ($\ref{eq.tnd}$) makes sense without the necessity of the
equivalence relation given by $\pmod{1}$.

In the random case, because the probability of trajectories which initially
rotate faster than the average rotation number is positive, the equality
between the rotation number of the discrete sampled system and the
continuous system only happens when we take the limit of the period $T$
going to zero, {\it i.e.} in the random case it does not exist a Nyquist's
rate ({\em cf.} Example 4).

\subsection{The Real Noise Case}

Consider a linear random equation on ${\R}^{2}$: 
\[
\dot{x}_{t}=A(\theta _{t}(\omega ))x_{t}\ ; 
\]
its fundamental solution $\varphi (t,\omega )$ is a continuous linear
cocycle. The continuous angular coordinate $\alpha _{t}$ of the solution $%
x_{t}$ satisfies the random equation: 
\[
\dot{\alpha}_{t}=\langle v_{t},A(\theta _{t}(\omega )s_{t}\rangle \ , 
\]
where $s_{t}$ is the radial projection of the solution $x_{t}$ on the
unitary circle $S^{1}$ and $v_{t}$ is such that $(s_{t},v_{t})$ is an
orthonormal pair with anti--clockwise orientation. By the ergodic theorem
the rotation number $\rho (\varphi )$ of this system satisfies: 
\[
\rho (\varphi )=\lim_{t\rightarrow \infty }\ \frac{1}{t}\ \int_{0}^{t}\ \dot{%
\alpha}_{r}\ dr={\rm I\!E}\left[ \int_{S^{1}}\langle v,A(\omega )s\rangle \
d\nu _{\omega }(s)\right] \ \ \ \ \ a.s., 
\]
where $\mu =\nu _{\omega }(ds){\rm I\!P}(d\omega )$ is an invariant
probability measure on $S^{1}\times \Omega $. So, formula ($\ref{rot.cocycle}
$), hence Theorem \ref{7-sampling thm}, holds in this case.

\subsection{The Stochastic Case}

Consider the following stochastic linear system in ${\R}^{2}$: 
\begin{equation}
dx_{t}=Ax_{t}\ dt+\sum_{i=1}^{m}B^{i}x_{t}\ \circ dW_{t}^{i}\ \ \ ,
\label{eq.linearsyst}
\end{equation}
where $A,B^{1},\ldots ,B^{m}$ are $2\times 2$--matrices, $(W_{t}^{1},\ldots
,W_{t}^{m})_{t\geq 0}$ is a Brownian motion in ${\R}^{m}$ with respect
to the probability space $(\Omega ,{\cal F},{\rm I\!P})$ with its natural
filtration $\{{\cal F}_{t}\}_{t\geq 0}$, and the integral is taken in the
Stratonovich sense. This stochastic system generates a white noise linear
cocycle $\varphi (t,\omega )$ (see {\it e.g.} Arnold \cite{Arnold}). The
continuous angular coordinate $\alpha _{t}\in {\R}$ of a solution $x_{t}$
with initial condition $x_{0}\in {\R}^{2}\setminus \{0\}$ is given by
the It\^{o} equation: 
\[
d\alpha _{t}=f(s_{t})\ dt+\sum_{i=1}^{m}<B^{i}s_{t},v_{t}>\ dW_{t}^{i}\ . 
\]
where $f:S^{1}\rightarrow {\R}$ is given by: 
\begin{equation}
f(s)=<As,v>+\sum_{i=1}^{m}\left( \frac{1}{2}%
<(B^{i})^{2}s,v>-<B^{i}s,s><B^{i}s,v>\right) \ ,  \label{eq.function}
\end{equation}
with $v$ such that $(s,v)$ is an orthonormal pair with anti--clockwise
orientation, see Ruffino \cite{Ruf}. In the stochastic case the invariant
measure $\mu $ on $\Omega \times S^{1}$ factorizes trivially a.s. $\mu ={\rm %
I\!P}\otimes \nu $, where $\nu $ is a stationary probability measure on $%
S^{1}$ for the Markov processes $s_{t}$. Hence, since the average of the
martingale part vanishes a.s., again, by the ergodic theorem the rotation
number satisfies: 
\[
\rho (\varphi )=\lim_{t\rightarrow \infty }\ \frac{1}{t}\ \int_{0}^{t}\
f(s_{r})\ dr=\int_{S^{1}}f(s)\ d\nu (s)\ \ \ \ \ \ \ \ a.s.. 
\]
Hence the condition given by equation ($\ref{rot.cocycle}$) is satisfied
considering the infinitesimal generator of the processes $\alpha _{t}$, so
Theorem \ref{7-sampling thm} also holds in this case.

\vspace{5mm}

We finish with an example which illustrates the fact that in the random case
it does not exist the Nyquist's rate. \vspace{2mm}\newline
{\sc Example 4:} Consider the following stochastic linear system: 
\[
dx_{t}=\left( 
\begin{array}{rr}
0 & -1 \\ 
1 & 0
\end{array}
\right) x_{t}\ dt+\left( 
\begin{array}{rr}
0 & -1 \\ 
1 & 0
\end{array}
\right) x_{t}\ \circ dW_{t}\ , 
\]
where $x_{t}\in {\R}^{2}$ and $(W_{t})_{t\geq 0}$ is a linear Brownian
motion. The fundamental solution of this equation $(\varphi _{t})_{t\geq 0}$
is the random linear rotation: 
\[
\varphi _{t}=\left( 
\begin{array}{rr}
\cos (t+W_{t}) & -\sin (t+W_{t}) \\ 
\sin (t+W_{t}) & \cos (t+W_{t})
\end{array}
\right) \ , 
\]
hence the continuous rotation $\alpha _{t}$ is independent of the starting
point $s\in S^{1}$ and is given by: 
\[
\alpha _{t}=t+W_{t}\ ; 
\]
therefore the rotation number 
\[
\rho =\lim_{t\rightarrow \infty }\ \ \frac{t+W_{t}}{t}\ =\ 1\ \ \ \ \ \ \
a.s. 
\]
The unique invariant probability $\nu $ in $S^{1}$ is the normalized
Lebesgue measure, besides, formula ($\ref{eq.function}$) gives $f(s)\equiv 1$%
, which confirms $\rho (\varphi )=1$ $\ a.s.$. For a discretization with
time interval $T>0$: 
\[
\beta _{1}^{T}(\omega ,s)\ =\ T+W_{T}\ \ \ \pmod{1} 
\]
with $\beta _{1}^{T}(\omega ,s)\in (-1/2,\ 1/2]$, and it does not depend on $%
s\in S^{1}$. The distribution of the random variable $\beta _{1}^{T}(\omega
) $ corresponds to the heat kernel $P_{T}(T,\cdot )$ in the circle $S^{1}$.
In the figure 1 the graphics with continuous curves show its distribution
for a sequence of decreasing values of $T$. Since the distribution of the
Brownian motion on $S^{1}$ is the canonical projection of the distribution
of the Brownian motion on ${\R}$ (the universal covering space of $S^{1}$%
) each continuous curve in the graphic are obtained by adding up the
projections of the distributions on unitary translated intervals of ${\R}
$ , represented by the dashed curves. Since in this case the rotation number 
$\rho (F_{T},F_{T}\circ \theta ,\ldots )={\rm I\!E}[\beta _{1}^{T}]\ $,
figure 1 makes clear the fact that 
\[
\frac{1}{T}\ \ \rho (F_{T},F_{T}\circ \theta ,\ldots )\ <\ 1\ \ \ a.s.\ \ , 
\]
for all $0<T<1/2$, so we do need the limit.

\begin{tabular}{cc}

\includegraphics[scale=0.25]{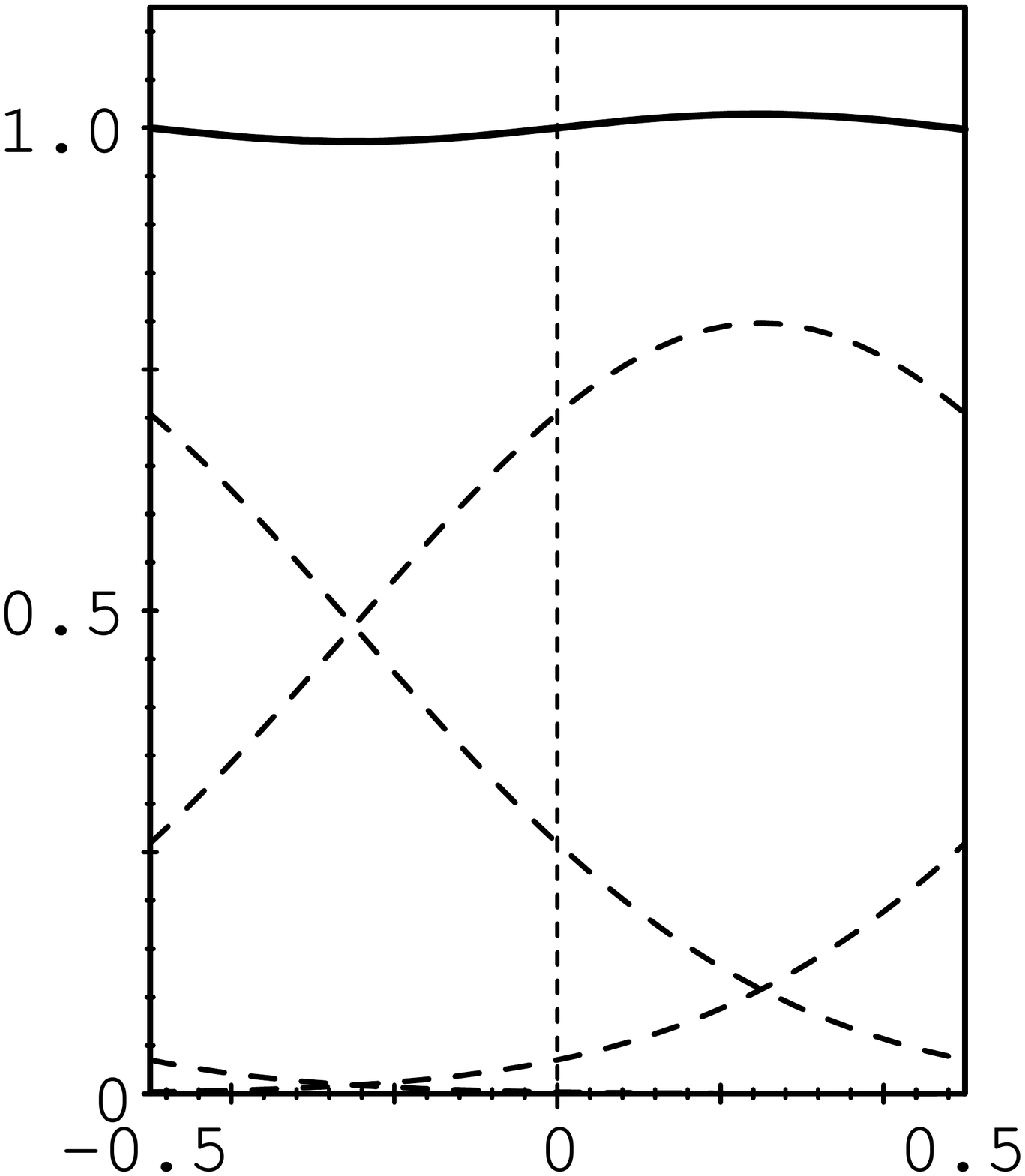}
 & \includegraphics[scale=0.25]{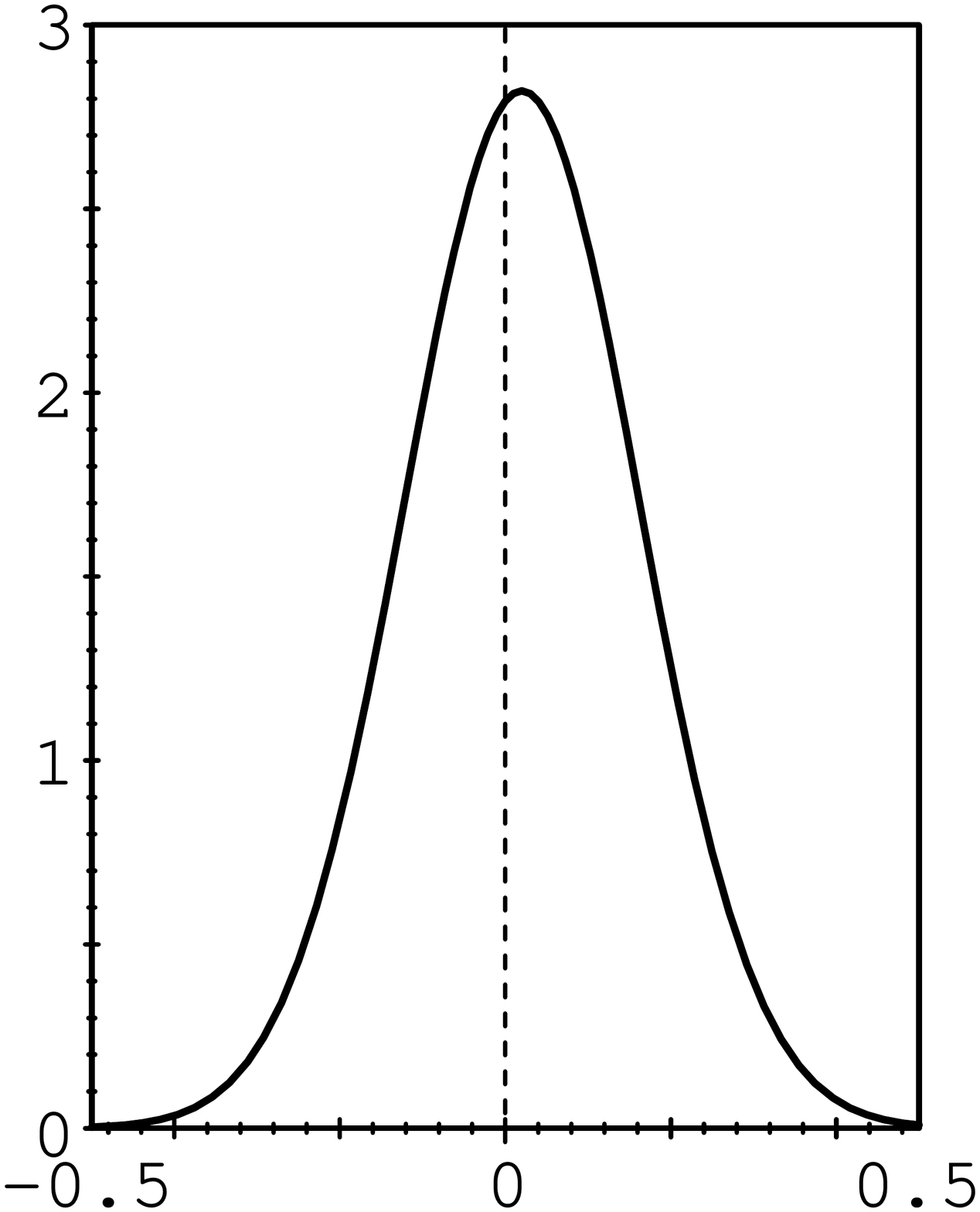}  \\
a) & b)  \\
${\rm I\!E}[\beta _{1}^{T}]\simeq 0.0259$ & ${\rm I\!E}[\beta _{1}^{T}]\simeq
0.0197$ \\
T=0.25 & T=0.02\\
\end{tabular}


\begin{center}
{\footnotesize Figure 1: Distribution of $\beta_1^T $ in the interval $%
(-1/2,\ 1/2]$ for decreasing values of $T$. }
\end{center}

\pagebreak

\begin{center}

\vspace{1cm}

 {\Large {\bf Erratum to: A sampling theorem for rotation numbers of 
linear processes in $\R^2$\\[3mm]
}}

\end{center}

\vspace{0.3cm}

\begin{center}
  Paulo R. Ruffino \footnote{E-mail:
ruffino@ime.unicamp.br.}

\vspace{0.2cm}

\textit{Departamento de Matem\'{a}tica, Universidade Estadual de Campinas, \\
13.083-859- Campinas - SP, Brazil.}

\end{center}

\vspace{20mm}

\noindent {\bf Erratum to:  Random Oper. Stochastic Equations 8 (2000), no. 2, 
175--188.}

\noindent {(To appear in Random Oper. Stochastic Equations 22 (2014), no. 
4.}

\setcounter{section}{0}

\section{Introduction}

The purpose of these notes is to correct the arguments in the proof of the 
random sampling theorem (Theorem 5.1) in the stochastic case 
(Subsection 5.3). The problem in the original proof, where we have claimed the 
convergence of the limit using just Lebesgue's convergence theorem, is that, in 
fact, there is an extra summand (see below) which is 
not obvious that it goes to zero. Precisely, in lines 13-14 from bottom to top 
in page 184, we have that 
\begin{eqnarray*} 
 &  \displaystyle \lim_{T\rightarrow 0} \  \frac{1}{T}\ \E \left[ \int_{S^1} 
\beta_1^T(\omega , s_0)\ \mathrm{d}\nu (s) 
\right] = & \\ 
&& \\
 & \displaystyle \lim_{T \rightarrow 0}\ \ \frac{1}{T} \ \E 
\left[ \int_{S^1} \left( \alpha_T(\omega) -\alpha_0 \right) \ 1_{T\leq \sigma } 
\ \mathrm{d} \nu (s)\right] &\\
&& \\
 & + \displaystyle \lim_{T \rightarrow 0}\ \ \frac{1}{T} \ \E 
\left[ \int_{S^1} \left( \alpha_T(\omega) -\alpha_0 \right) \ 1_{T > \sigma } 
\ \mathrm{d} \nu (s)\right].
\end{eqnarray*}
Where the last summand in the right hand side was na\"ively considered as zero, 
only based on the fact that the numerator goes to zero. Recently we have 
realized that depending on the (nongaussian) noise, this term may not vanish.
Hence an extra care has do be done. Our intention here is to complete the 
argument, proving that, in fact, for stochastic systems  (Gaussian noise) this 
last term does converge to zero, establishing, then the sampling 
theorem for this case.

\bigskip

\noindent {\bf Proof:}  Recalling our notations, for a given  stochastic 
linear equation in $\R^2$, the 
continuous angular coordinate $\alpha_t$ of this system satisfies
\begin{equation} \label{eq: para alpha}
d\alpha_t  =  f(s_t) \mathrm{d}t +  + \sum_{i=1}^m  <B^i s_t, v_t>\ 
\mathrm{d} W_t^i.
\end{equation}
where $s_t= \pi (\alpha_t)$, $v_t$ is 
orthonormal to $s_t$ with positive orientation and $f: S^1 \rightarrow \R$ is 
given by 
\[
f(s)= <As, v> + \sum_{i=1}^m \left( \frac{1}{2} <(B^i)^2s, v> - <B^is, s> 
<B^is, v> \right).
\]

 For an initial condition $s_0 \in (-1/2, 1/2] \sim S^1$ and all 
$T>0$, we have that,
\begin{equation} \label{eq: beta}
 \beta^T_1 (\omega, s_0) = \left( \alpha_T (\omega) - \alpha_0 (\omega) \right) 
+ N(\omega)
\end{equation}
for an integrable integer variable $N(\omega)$. According to our 
construction, $N(\omega)$ only depends on the trajectory of $\alpha_0=0$, or of 
$e_1\in S^1$. It measures how many times this trajectory crosses its antipode 
$(-e_1)$ in the 
anti-clockwise direction during the interval $[0,T]$.  When $T> \tau 
(\omega)$, we have that $N(\omega) \neq 0$.

Hence, the proof is completed if we control the expectation $\E 
[|N|]$ and show 
that it goes to zero faster than $T$. We use  boundedness on the
distribution of $\alpha_t$ with initial condition $\alpha_0=0$. Let $p(t, x,y)$ 
be the
density of the transition probability measure associated to the
non-degenerate diffusions given by Equation 
(\ref{eq: para alpha}). 
Then, there exists a constant $M>0 $ such 
that,

\[
\frac{1}{M \sqrt{t}} e^{-M\frac{(x-y)^2}{t}} \leq p (t, x,y) \leq \frac{M}{ 
\sqrt{t}} e^{-\frac{(x-y)^2}{Mt}}.
\]
See Kusuoka and Stroock \cite{Kusuoka and Stroock, Kusuoka and Stroock-1988}. 
Let $N^+= \max \{N, 0\}$ and $N^- =  \max \{-N, 0 \} $, such that 
$N=N^+ - N^-$. Hence, for the positive 
part $N^+$
\begin{eqnarray*}
 \E [ N^+ ] & \leq & M \int_{\alpha +1}^{\infty} \left( \lfloor x-(\alpha+1) 
\rfloor + 1 \right) \frac{1}{\sqrt{T}} \exp\left\{- \frac{ 
(x-q)^2}{M T}\right\} \ dx \\
  && \\
  & \leq & M \int_{\alpha +1}^{\infty} \left( x- \alpha \right) 
\frac{1}{\sqrt{T}} \exp\left\{- \frac{ 
(x-q)^2}{M T}\right\} \ dx.
\end{eqnarray*}
And for the negative part:
\begin{eqnarray*}
 \E [ N^- ] & \leq & M \int_{-\infty}^{\alpha} \left( \lfloor \alpha -x 
\rfloor + 1 \right) \frac{1}{\sqrt{T}} \exp\left\{ - \frac{ 
(x-q)^2}{M T}\right\} \ dx \\
  && \\
  & \leq & M \int_{-\infty}^{\alpha} \left(  \alpha -x + 1 \right) 
\frac{1}{\sqrt{T}} \exp\left\{ - \frac{ 
(x-q)^2}{M T}\right\} \ dx.
\end{eqnarray*}
Changing variables, for $T \in (0,1)$ we have that
\begin{eqnarray*}
 \E [ N^+ ] & \leq & M^{\frac{3}{2}} \int_{\frac{\alpha -q 
+1}{\sqrt{M T}}}^{\infty} \left( \sqrt{M} u + q - \alpha  \right) 
 \exp\left\{- u^2 \right\} \ du
\end{eqnarray*}
and
\begin{eqnarray*}
 \E [ N^- ] & \leq & M^{\frac{3}{2}} \int_{\frac{ q -\alpha}{\sqrt{M 
T}}}^{\infty} \left( \sqrt{M} u + \alpha -q + 1  \right) 
 \exp\left\{- u^2 \right\} \ du
\end{eqnarray*}

Hence $\E [ N^+ ]$ and $\E [ N^- ]$ goes to zero when $T$ goes to zero. 
Moreover, by standard calculus argument, 
using that $\lim_{z\rightarrow 0} 
\exp\left\{ {-\frac{1}{z}} \right\} z^{\beta}=0$ for any exponent $\beta \in 
\R$, then finally we 
get that
\[
  \lim_{T \searrow 0} \frac{\E [N^+ ]}{T} =  \lim_{T 
\searrow 0} \frac{\E [N^- ]}{T} = 0.
\]

\eop

\bigskip

\noindent {\bf Acknowledgments:} The author would like to thank discussion with 
Dr. Christian Rodrigues, who pointed out that for the Gaussian stochastic case, 
 Theorem 5.1 needs extra arguments. The 
author is partially supported by FAPESP 11/50151-0 and
12/03992-1.

\end{document}